\documentclass{amsart}

\usepackage{amssymb}
\usepackage{stmaryrd}
\usepackage{enumerate}
\usepackage{hyperref}
\usepackage{fancybox}
\usepackage{mathrsfs}
\usepackage[matrix,arrow,graph]{xy}

\theoremstyle{plain}      \newtheorem*{theorem}{Theorem}
\theoremstyle{plain}      
\theoremstyle{plain}      
\theoremstyle{plain}      
\theoremstyle{definition} 
\theoremstyle{definition} 
\theoremstyle{remark}     
\theoremstyle{remark}     

\renewcommand{\mathcal}{\mathscr}
\newcommand{\A}{\ensuremath{\mathcal{A}}}
\newcommand{\B}{\ensuremath{\mathcal{B}}}
\newcommand{\C}{\ensuremath{\mathcal{C}}}
\newcommand{\D}{\ensuremath{\mathcal{D}}}

\newcommand{\V}{\ensuremath{\mathcal{V}}}

\newcommand{\Prof}{\ensuremath{\mathbf{Prof}}}

\newcommand{\x}{\ensuremath{\times}}
\renewcommand{\o}{\ensuremath{\circ}}
\newcommand{\bs}{\ensuremath{\backslash}}
\newcommand{\ox}{\ensuremath{\otimes}}
\newcommand{\op}{\ensuremath{\mathrm{op}}}
\newcommand{\ob}{\ensuremath{\mathrm{ob}}}

\newcommand{\ra}{\ensuremath{\xymatrix@=4ex@1{\ar[r]&}}}
\newcommand{\Ra}{\ensuremath{\xymatrix@=4ex@1{\ar@{=>}[r]&}}}
\renewcommand{\mapsto}{\ensuremath{\xymatrix@=4ex@1{\ar@{|->}[r]&}}}

\begin{document}

\title{Biclosed bicategories: localisation of convolution}
\author{Brian J. Day}
\address{Centre of Australian Category Theory, Macquarie University, NSW,
2109, Australia}
\date{\today}

\begin{abstract}
We give a summary (without proofs) of the main results in the author's
thesis entitled ``Construction of biclosed categories'' (University of New
South Wales, Australia, 1970). This summary is reprinted directly from
Report~81-0030 of the School of Mathematics and Physics, Macquarie
University, April 1981. In particular, it gives sufficient conditions for
existence of an extension of a (pro)monoidal category structure along a
given dense functor to a cocomplete category. The two basic procedures used in
the proof turn out to be special cases of the final result, the two
respective dense functors then being the Yoneda embedding followed by a
localisation. The final result has a standard universal property based on
left Kan extension of (pro)monoidal functors along the given dense functor,
however this property is not stated explicitly here.
\end{abstract}

\maketitle

\section{Introduction}

The aim of this article is to record, for future reference, some of the
elementary formulas arising in the theory of biclosed bicategories. These
formulas are given by two outstanding building blocks in the theory, namely
convolution~\cite{2} and reflection~\cite{3}. When combined, these processes
lead to an extension theorem for probicategories (which is the name we give
to the collection of structure functors of the basic biclosed bicategories
of functors (or ``presheaves'') into the ground category).

All categorical algebra shall be relative to a fixed complete and cocomplete
symmetric monoidal closed ground category $\V = (\V,\ox,I,[-,-],\ldots)$.

\section{Probicategories and convolution}

A \emph{probicategory} structure on a family $\{\A_{xy} : x,y \in \ob(\A)\}$
of small categories is essentially a biclosed bicategory structure on the
collection $\{[\A_{xy},\V] : x,y \in \ob(\A)\}$, each of the composition
functors
\[
\o : [\A_{yz},\V] \ox [\A_{xy},\V] \ra [\A_{xz},\V]
\]
being determined, to within isomorphism, by a \emph{structure} functor:
\[
P_{xyz}:\A^\op_{yz} \ox \A^\op_{xy} \ox \A_{xz} \ra \V.
\]
A bicategory is a particular instance of a probicategory for which there
exist functors
\[
\o:\A_{yz} \ox \A_{xy} \ra \A_{xz}
\]
such that
\[
P_{xyz}(A,B,C) \cong \A_{xz}(A \o B, C).
\]
A biclosed probicategory is a probicategory for which there exist functors
\begin{align*}
-/- &: \A_{xz} \ox \A^\op_{xy} \ra \A_{yz} \\
- \bs- &: \A^\op_{yz} \ox \A_{xz} \ra \A_{xy}
\end{align*}
such that
\begin{align*}
P_{xyz}(A,B,C) &\cong \A_{xy}(B,A \bs C) \\
               &\cong \A_{yz}(A,C / B)
\end{align*}
A manifold probicategory is an indexing of $\Prof$. Thus a set $\ob(\A)$ is
given and, for each $x \in \ob(\A)$, a category $\A_x$. We take $\A_{xy} =
\A^\op_x \ox \A_y$ and
\[
P_{xyz}((A,A'),(B,B'),(C,C')) = \A_x(C,B) \ox \A_y(B',A) \ox \A_z(A',C').
\]
A probicategory may (or may not) have an \emph{identity} $J$; we shall here
include the identity conditions throughout.

One can always convolve a probicategory $\A = (\A,P,J,\ldots)$ with a
cocomplete monoidal category $\B = (\B,\ox,I,\ldots)$ for which $- \ox B$
and $B \ox -$ both preserve colimits for all objects $B \in \B$. The
bicategory composition and identities on the convolution $[\A,\B]$ are given
by:
\begin{align*}
F \o G &= \int^{AA'} P(A,A',-) * (FA \ox GA') \\
I &= J * I.
\end{align*}

\section{Reflection theorem}

A morphism $\psi = (\psi,\tilde{\psi},\psi^\o):\B \ra \C$ of bicategories
will be called \emph{strong} if each of the comparison maps
$\tilde{\psi}:\psi B \o \psi B' \ra \psi(B \o B')$ and $\psi^\o:I \ra \psi
I$ is an isomorphism.

\begin{theorem}
Let $\B = (\B_{xy}, \o, I,/,\bs,\ldots)$ be a biclosed bicategory and let
$\theta:\C \ra \B$ be a family of full embeddings $\theta_{xy}:\C_{xy} \ra
\B_{xy}$ with left adjoints $\psi = \{\psi_{xy}\}$ (we omit $\theta$ from
the notation, and denote the unit of the adjunction by $\eta:1 \ra \psi:\B
\ra \B$). Further, let $\A_{xy} \subset \B_{xy}$ be a locally strongly
generating class of 1-cells in $\B$ and let $\D_{xy} \subset \C_{xy}$ be a
locally strongly cogenerating class of 1-cells in $\C$. Then, in order that
there should exist a biclosed bicategory structure on $\C$ for which
$\psi:\B \ra \C$ admits enrichment to a strong map of biclosed bicategories,
it is
sufficient that \emph{one} of the following pairs of morphisms be a pair of
isomorphisms for all appropriately indexed 1-cells $A \in \A$, $B,B' \in \B$,
$C \in \C$, and $D \in \D$:
\[
\begin{array}{lll}
\text{1.} & \text{a)} & \eta: C/B \ra \psi(C/B) \\
         & \text{b)} & \eta: B \bs C \ra \psi(B \bs C)
\medskip\\
\text{2.} & \text{a)} & \eta: D/A \ra \psi(D/A) \\
         & \text{b)} & \eta: A \bs D \ra \psi(A \bs D) 
\medskip\\
\text{3.} & \text{a)} & \eta \bs 1: \psi B \bs C \ra B \bs C \\
         & \text{b)} & 1/\eta: C/\psi B \ra C / B 
\medskip\\
\text{4.} & \text{a)} & \psi(\eta \o 1): \psi(B \o B') \ra \psi(\psi B \o B')\\
         & \text{b)} & \psi(1 \o \eta): \psi(B' \o B) \ra \psi(B' \o \psi B)
\medskip\\
\text{5.} & \text{a)} & \psi(\eta \o 1): \psi(B \o A) \ra \psi(\psi B \o A)\\
         & \text{b)} &  \psi(1 \o \eta): \psi(A \o B) \ra \psi(A \o \psi B)
\medskip\\
\text{6.} & & \psi(\eta \o \eta): \psi(B \o B') \ra \psi(\psi B \o \psi B').
\end{array}
\]
\end{theorem}

A special application of the reflection theorem is to the localisation of a
probicategory $\A$ (cf.~\cite{4}). Let $\Sigma$ be a set of 2-cells in $\A$
and form the locally Cauchy-dense map:
\[
\Pi:\A \ra \A(\Sigma^{-1}).
\]
Then the conditions of the reflection theorem applied to the full reflective
embedding:
\[
\psi \dashv [\Pi,1]:[\A(\Sigma^{-1}),\V] \ra [\A,\V]
\]
are conditions for $\A(\Sigma^{-1})$ to carry a probicategory structure
such that $\Pi$ becomes a ``map of probicategories''.

\section{Extension theorem}

The following result is obtained by locally completing the family $\C$ with
respect to a suitable change of $\V$-universe.

\begin{theorem}
Let $\A = (\A,P,J,\ldots)$ be a probicategory and let $N_{xy}:\A^\op_{xy}
\ra \C_{xy}$ be an $\ob(\A) \x \ob(\A)$-indexed family of dense functors. Then
$\C$ admits a biclosed bicategory structure, extending the structure of
$\A$, if the following indexed colimits and limits exist in $\C$ for all
appropriately indexed 1-cells $A,A' \in \A$ and $C,C' \in \C$:
\begin{align*}
Q_{xyz}(A,A') &= P_{xyz}(A,A',X)*N_{xz} X, \\
I &= J_{xx} X * N_{xx} X, \\
C \o C' &= (\C_{yz} (N_{yz} X,C) \ox \C_{xy}(N_{xy} X',C') * Q_{xyz}(X,X'),\\
H_{yz}(A,C) &= \C_{xz}(Q_{xyz}(X,A),C) * N_{yz} X, \\
K_{xy}(A,C) &= \C_{xz}(Q_{xyz}(A,X),C) * N_{xy} X, \\
C/C' &= \{\C_{xy}(N_{xy} X,C'),H_{yz}(X,C)\}, \\
C \bs C' &= \{\C_{yz}(N_{yz} X,C),K_{xy}(X,C')\},
\end{align*}
and the following induced natural transformations
\begin{align*}
\xymatrix{\C_{xz}(Q_{xyz}(A',A),C) \ar[r]^-\cong &
    \C_{yz}(N_{yz}A',H_{yz}(A,C))} \\
\xymatrix{\C_{xz}(Q_{xyz}(A,A'),C) \ar[r]^-\cong &
    \C_{xy}(N_{xy}A',K_{xy}(A,C))}
\end{align*}
are isomorphisms (where
\begin{align*}
H_{yz}(-,-) &: \A_{xy} \ox \C^\op_{xz} \ra \C_{yz} \\
K_{xy}(-,-) &: \A^\op_{yz} \ox \C_{xz} \ra \C_{xy}, \text{ and} \\
-/- &: \C_{xz} \ox \C^\op_{xy} \ra \C_{yz} \\
- \bs - &: \C^\op_{yz} \ox \C_{xz} \ra \C_{xy}).
\end{align*}
\end{theorem}




\begin{thebibliography}{10}

\bibitem{B}
J. B\'enabou, Introduction to bicategories, Lecture Notes in Math. 47
Springer (1967) pp 1--77.

\bibitem{2}
B. Day, On closed categories of functors, Reports of the Midwest Category
Seminar IV, Lecture Notes in Mathematics, Vol. 137 (Springer 1970), 1--38.

\bibitem{3}
B. Day, A reflection theorem for closed categories, J. Pure Appl. Algebra 2
(1972), 1--11.

\bibitem{4}
B. Day, Note on monoidal localisation, Bull. Austral. Math. Soc. 8 (1973),
1--16.

\end{thebibliography}
\end{document}